\pdfoutput=1
\RequirePackage{ifpdf}
\ifpdf 
\documentclass[pdftex]{sigma}
\else
\documentclass{sigma}
\fi

\numberwithin{equation}{section}
\newtheorem{Theorem}{Theorem}[section]
\newtheorem{Proposition}[Theorem]{Proposition}
 { \theoremstyle{definition}
\newtheorem{Definition}[Theorem]{Definition}
\newtheorem{Remark}[Theorem]{Remark} }

\providecommand{\co}{\colon\,}

\begin{document}

\allowdisplaybreaks

\renewcommand{\thefootnote}{$\star$}

\renewcommand{\PaperNumber}{071}

\FirstPageHeading

\ShortArticleName{Levi-Civita's Theorem for Noncommutative Tori}

\ArticleName{Levi-Civita's Theorem for Noncommutative Tori\footnote{This paper is a~contribution to the
Special Issue on Noncommutative Geometry and Quantum Groups in honor of Marc A.~Rief\/fel.
The full collection is available at \href{http://www.emis.de/journals/SIGMA/Rieffel.html}
{http://www.emis.de/journals/SIGMA/Rief\/fel.html}}}

\Author{Jonathan ROSENBERG}

\AuthorNameForHeading{J.~Rosenberg}
\Address{Department of Mathematics, University of Maryland, College Park, MD 20742, USA}
\Email{\href{mailto:jmr@math.umd.edu}{jmr@math.umd.edu}}
\URLaddress{\url{http://www.math.umd.edu/~jmr/}}

\ArticleDates{Received July 26, 2013, in f\/inal form November 19, 2013; Published online November 21, 2013}

\Abstract{We show how to def\/ine Riemannian metrics and connections on a~noncommutative torus in such
a~way that an analogue of Levi-Civita's theorem on the existence and uniqueness of a~Riemannian connection
holds.
The major novelty is that we need to use two dif\/ferent notions of noncommutative vector f\/ield.
Levi-Civita's theorem makes it possible to def\/ine Riemannian curvature using the usual formulas.}


\Classification{46L87; 58B34; 46L08; 46L08}

\begin{flushright}
\it To Marc Rieffel, with admiration and appreciation
\end{flushright}

\renewcommand{\thefootnote}{\arabic{footnote}} \setcounter{footnote}{0}

\pdfbookmark[1]{Introduction}{intro}
\section*{Introduction}

In his lecture series at the Focus Program on Noncommutative Geometry and Quantum Groups at the Fields
Institute in June, 2013, Masoud Khalkhali gave a~very beautiful description of recent work by Connes and
Moscovici~\cite{2011arXiv1110.3500C} (building on earlier work of Connes and Tretkof\/f~\cite{MR2907006})
and by Fathizadeh and Khalkhali~\cite{MR2956317,2011arXiv1110.3511F,2013arXiv1301.6135F} on a~calculation
of what one can call ``scalar curvature'' for metrics on noncommutative tori obtained by (noncommutative)
conformal deformation of a~f\/lat metric.
At the same time, Khalkhali explained that def\/ining curvature in terms of the spectral geometry of the
Laplacian is basically forced on us by a~lack in the noncommutative setting of the standard machinery of
Riemannian geometry, whereby one would def\/ine curvature using derivatives of the Levi-Civita connection.
In the same Focus Program, Marc Rief\/fel in his lecture gave a~def\/inition of Riemannian metric in the
noncommutative setting, albeit only for f\/inite-dimensional algebras (which, roughly speaking, correspond
to zero-dimensional manifolds).
The purpose of this paper is to show that one can give a~very natural and quite general def\/inition of
Riemannian metrics on a~noncommutative torus, without assuming \emph{a priori} that the metric is
a~conformal deformation of a~f\/lat metric, and that for this def\/inition one can prove an analogue of
Levi-Civita's theorem (\cite{zbMATH02613477}; for a~modern formulation and proof see for example~\cite[Chapter~2, \S~3]{MR1138207}) on the existence and uniqueness of a~torsion-free connection compatible with the metric.
For our notion of Riemannian metric we also obtain a~notion of Riemannian curvature, which we compute
explicitly in the two-dimensional case.

Admittedly our def\/inition of Riemannian metric still has certain drawbacks.
It would be nice to be able to prove a~uniformization theorem for two-dimensional noncommutative tori,
stating that in some sense all Riemannian metrics are equivalent to conformal deformations of a~standard
f\/lat metric, and are thus of the form studied by Connes--Moscovici and by Fathizadeh--Khalkhali.
The problem, however, is that generic smooth two-dimensional noncommutative tori are quite rigid~-- their
dif\/feomorphism groups are not much bigger than the group of smooth inner
automorphisms~\cite{MR859436}~-- and this makes it hard to see how a~uniformization theorem could be true
without using a~very dif\/ferent def\/inition of Riemannian metric.

\vspace{-1mm}

\section{Riemannian metrics and connections}

While it would be desirable to have a~theory of Riemannian metrics and connections on arbitrary
``noncommutative manifolds,'' there is a~problem in general in understanding what a~tangent vector or
vector f\/ield should be.
However, tori are parallelizable, so on a~torus, a~vector f\/ield is simply a~linear combination (with the
coef\/f\/icients being arbitrary functions) of the (commuting) coordinate vector f\/ields~$\partial_j$.
This def\/inition carries over without dif\/f\/iculty to a~noncommutative torus, though as we will see,
there are also other ways of def\/ining vector f\/ields.
We begin with basic def\/initions and notation.
\begin{Definition}
Let $\Theta$ be a~skew-symmetric $n\times n$ matrix with entries in $\mathbb{R}$.
(When $n=2$, we write $\Theta=
\begin{pmatrix}
0&\theta
\\
-\theta&0
\end{pmatrix}
$ for $\theta\in \mathbb{R}$.) The \emph{noncommutative torus} of dimension $n$ with noncommutativity
parameter $\Theta$ is the universal {$C^*$-algebra} $A_\Theta$ on $n$ unitaries $U_1,\dots,U_n$ with
commutation relations $U_jU_k = e^{2\pi i \Theta_{jk}}U_kU_j$.
(When $n=2$, there are only two generators $U_1$ and $U_2$ and there is only one relation, $U_1U_2 =
e^{2\pi i \theta}U_2U_1$, and we write $A_\theta$ for $A_\Theta$.) We will sometimes think of $A_\Theta$ as
the algebra of ``functions'' on a~``noncommutative manifold'' $\mathbb{T}^n_\Theta$.
The algebra $A_\Theta$ carries an ergodic action of the $n$-torus $\mathbb{T}^n$ via $t\cdot U_j = t_jU_j$,
$t\in \mathbb{T}^n$.
The inf\/initesimal generators of this action are the (unbounded) $*$-derivations $\partial_j$ with
$\partial_j(U_k) = \delta_{jk}2\pi i U_k$.
The $C^\infty$ vectors for this action constitute the algebra $A^\infty_\Theta$ called the \emph{smooth
noncommutative torus}.
This algebra can be identif\/ied with the (noncommutative) rapidly decreasing Fourier series
\begin{gather*}
\left\{\sum_{m_1,\dots,m_n}c_{m_1,\dots,m_n}U_1^{m_1}\cdots U_n^{m_n}\Big|\{c_{m_1,\dots,m_n}\}
\ \text{rapidly decreasing}\right\}.
\end{gather*}
When $\Theta=0$ (the commutative case), this is $C^\infty(\mathbb{T}^n)$.
In general this is always indistinguishable from $C^\infty(\mathbb{T}^n)$ as a~Fr\'echet space, even though
the algebra structures are dif\/ferent.
\end{Definition}

Now we get to the point where the noncommutative theory diverges from the usual geometry of manifolds.
On a~smooth (ordinary) closed manifold, there are three equivalent ways of def\/ining vector f\/ields: as
inf\/initesimal generators of one-parameter groups of dif\/feomorphisms (i.e., f\/lows), or as smooth
sections of the tangent bundle, or as homogeneous linear f\/irst-order dif\/ferential operators
annihilating the constants.
The problem is that in the noncommutative case, these def\/initions do not all agree, and so we need more
than one notion.
\begin{Definition}
\label{def:VF}
The space $\mathcal{X}_\Theta$ of \emph{vector fields} on the smooth noncommutative torus is def\/ined to
be the free rank-$n$ left $A^\infty_\Theta$-module with basis $\partial_1,\dots,\partial_n$.
In other words, a~vector f\/ield is just a~formal linear combination of ``partial derivatives'' with
``function'' coef\/f\/icients in $A^\infty_\Theta$.
Vector f\/ields operate on $A^\infty_\Theta$ in the obvious way as linear ``f\/irst-order dif\/ferential
operators'' annihilating the ``constant functions'' $\lambda\cdot 1$, $\lambda\in \mathbb{C}$.
\end{Definition}

The problem with this def\/inition is that, unlike the standard commutative situation, an element of
$\mathcal{X}_\Theta$ is \emph{not} usually a~derivation, and commutator of such vector f\/ields is not
necessarily a~vector f\/ield, since
\begin{gather*}
[b\partial_j,c\partial_k]a=b\partial_j(c\partial_k a)-c\partial_k(b\partial_j a)
=b(\partial_j c)(\partial_k a)+bc(\partial_j\partial_k a)-c(\partial_k b)(\partial_j a)-cb(\partial_k\partial_j a),
\end{gather*}
so{\samepage
\begin{gather*}
[b\partial_j,c\partial_k]=b(\partial_j c)\partial_k-c(\partial_k b)\partial_j+[b,c]\partial_j\partial_k
\end{gather*}
and the second order term $[b,c]\partial_j\partial_k$ does not necessarily cancel out.}

\begin{Definition}
Accordingly, we introduce another linear space $\mathcal{D}_\Theta$, consisting of the $*$-deriva\-tions
$\delta\co A^\infty_\Theta\to A^\infty_\Theta$.
By~\cite[Corollary~5.3,~C2]{MR727402}, any such $\delta$ is automatically continuous in the Fr\'echet
topology.
It is clear that $\mathcal{D}_\Theta$ is a~Lie algebra under the commutator bracket (since the bracket of
derivations is a~derivation), so this remedies one of the defects of Def\/inition~\ref{def:VF}.
We can view $\mathcal{D}_\Theta$ as the Lie algebra of the inf\/inite dimensional group
$\operatorname{Dif\/f} (\mathbb{T}^n_\Theta) = \operatorname{Aut} (A^\infty_\Theta)$ of $*$-automorphisms
of $A^\infty_\Theta$.
Furthermore, by~\cite[Corollary 5.3, D2]{MR727402}, any $\delta \in \mathcal{D}_\Theta$ has a~unique
decomposition as $a_1\partial_1+\dots + a_n\partial_n+\delta_0$, where $\delta_0$ is approximately inner
and $a_1,\dots, a_n$ lie in the center of $A_\Theta^\infty$.
\end{Definition}
Now we need the following result:
\begin{Theorem}[Bratteli, Elliott, and Jorgensen]\label{thm:BEJ}
If $\Theta$ is ``generically transcendental'' $($in a~rather complicated sense made precise
in~{\rm \cite{MR727402}}, but satisfied for almost all skew-adjoint matrices$)$, then any $\delta \in
\mathcal{D}_\Theta$ has a~unique decomposition as $a_1\partial_1+\dots + a_n\partial_n+\delta_0$, where
$\delta_0$ is {\bfseries inner}, hence bounded in the {$C^*$-algebra} norm, and $a_1,\dots, a_n\in
\mathbb{C}$.
\end{Theorem}

\begin{proof}See~\cite[Remark~4.3]{MR727402}.
It also follows that $\delta$ is a~pregenerator of a~one-parameter subgroup of $\operatorname{Aut}
(A^\infty_\Theta)$ (in fact sometimes this is even true without genericity of $\Theta$,
see~\cite[Theo\-rem~5.4]{MR727402}).
\end{proof}

\begin{Definition}
\label{def:Rmetric}
A \emph{Riemannian metric} $g=\langle \cdot, \cdot\rangle$ on $A^\infty_\Theta$ is def\/ined to be
a~(positive) $A^\infty_\Theta$-valued inner product on vector f\/ields, or in other words, a~sesquilinear
map $\langle \cdot, \cdot\rangle\co \mathcal{X}_\Theta \times \mathcal{X}_\Theta \to A^\infty_\Theta$
satisfying the axioms of a~(pre-)Hilbert module:
\begin{enumerate}\itemsep=0pt
\item $\langle X + X', Y\rangle = \langle X, Y\rangle + \langle X', Y\rangle$ and $\langle aX, Y\rangle =
a\langle X, Y\rangle$ for $X,X',Y\in \mathcal{X}_\Theta$, $a\in A^\infty_\Theta$ (so $g$ is
$A^\infty_\Theta$-linear in the f\/irst variable); \item $\langle X, Y\rangle^* = \langle Y, X\rangle$
(hermitian symmetry)~-- together with (1), this implies $\langle X, aY\rangle = \langle X, Y\rangle a^*$;
\item $\langle X, X\rangle \ge 0$ in the sense of the {$C^*$-algebra} $A_\Theta$, with equality only if
$X=0$.
\end{enumerate}

Note that the metric is uniquely determined by the matrix $(g_{jk}) = \left(\left\langle
\partial_j,\partial_k \right\rangle\right)$ in $M_n(A^\infty_\Theta)$.
This matrix must be a~positive element of $M_n(A_\Theta)$ since
\begin{gather*}
0\le\left\langle\sum_j a_j\partial_j,\sum_k a_k\partial_k\right\rangle
=\sum_{j,k}a_j g_{jk}a_k^* \qquad \text{for any} \ \ a_j\in A^\infty_\Theta.
\end{gather*}

Actually, the axioms so far only correspond to a~hermitian metric on the complexif\/ied tangent bundle.
Recall that a~Riemannian metric must assign a~real-valued (i.e., self-adjoint) inner product to two real
vector f\/ields.
Since the ``real vector f\/ields'' are generated by the $\partial_j$, we need to add one additional
condition:
\begin{enumerate}\itemsep=0pt
\setcounter{enumi}{3} \item For each $j$, $k$, $\langle X_j, X_k\rangle$ is self-adjoint, and thus $\langle
X_j, X_k\rangle = \langle X_k, X_j\rangle$.
\end{enumerate}
\end{Definition}

\begin{Definition}
\label{def:conn}
A \emph{connection} on $A^\infty_\Theta$ is a~way of def\/ining covariant derivatives for vector f\/ields
satisfying the Leibniz rule
\begin{gather}
\label{eq:conn}
\nabla_X(aY)=(X\cdot a)Y+a\nabla_{X}Y.
\end{gather}
But here's the tricky aspect of this.
For~\eqref{eq:conn} even to make sense, we need to be able to multiply~$Y$ on the left by an element of the
algebra, so we want $Y\in \mathcal{X}_\Theta$.
On the other hand, applying~\eqref{eq:conn} to $\nabla_X\bigl((ab)Y\bigr)$ and comparing with the expansion
of $\nabla_X\bigl(a(bY)\bigr)$, we obtain
\begin{gather*}
(X\cdot(ab))Y=\bigl((X\cdot a)b+a(X\cdot b)\bigr)Y,
\end{gather*}
and since $Y$ is arbitrary, this forces $X$ to be a~derivation.
So we need $X\in \mathcal{D}_\Theta$ rather than $X\in \mathcal{X}_\Theta$.

In other words, a~connection is a~map $\nabla\co \mathcal{D}_\Theta \times \mathcal{X}_\Theta \to
\mathcal{X}_\Theta$, written $(X,Y)\mapsto \nabla_X Y$, satisfying the following axioms:
\begin{enumerate}\itemsep=0pt
\item $\nabla$ is linear in the f\/irst variable, so that $\nabla_{\lambda X}Y = \lambda\nabla_{X}Y$ and
$\nabla_{X+X'}Y = \nabla_{X}Y + \nabla_{X'}Y$, for $X,X'\in \mathcal{D}_\Theta$, $\lambda\in \mathbb{C}$;
\item $\nabla$ is $\mathbb{C}$-linear in the second variable, so $\nabla_X (Y+Y') = \nabla_X Y + \nabla_X
Y'$, $\nabla_X (\lambda Y) = \lambda \nabla_X Y$, for $X,Y,Y'\in \mathcal{X}_\Theta$, $\lambda\in
\mathbb{C}$; \item For $X,Y\in \mathcal{X}_\Theta$, $a\in A^\infty_\Theta$, $\nabla_X(aY) = (X\cdot a)Y +
a~\nabla_{X}Y$.
\end{enumerate}
Normally (i.e., in the classical case $\Theta=0$) the axioms for a~connection require that $\nabla$ be
$A^\infty_\Theta$-linear in the f\/irst variable, but this does not make sense in our context since
$\mathcal{D}_\Theta$ is not a~left $A^\infty_\Theta$-module.
However, let's assume that $\Theta$ is generic in the sense of Theorem~\ref{thm:BEJ}, so that any element
of $\mathcal{D}_\Theta$ dif\/fers from a~linear combination of $\partial_1,\dots, \partial_n$ by an inner
derivation.
We need an extra axiom to pin down the values of $\nabla_{\operatorname{ad} a}$, $a\in A^\infty_\Theta$.
We have no classical precedent for this since there are no inner derivations in the commutative case, but
from~\eqref{eq:conn} we obtain
\begin{gather*}
\nabla_{\operatorname{ad}a}(bY)=[a,b]Y+b\nabla_{\operatorname{ad}a}Y,
\qquad
\text{or}
\qquad
[\nabla_{\operatorname{ad}a},b]=[a,b].
\end{gather*}
The easiest way to satisfy this is to take $\nabla_{\operatorname{ad} a} =$ left multiplication by $a$.
However $\operatorname{ad} a$ only determines $a$ up to addition of a~constant, so we use the canonical
trace $\tau$ on $A_\Theta$ to normalize things.
Given the derivation $\operatorname{ad} a$, $a$ is unique subject to the condition that $\tau(a)=0$, and we
add as another axiom:
\begin{enumerate}\itemsep=0pt
\setcounter{enumi}{3} \item For any $a\in A^\infty_\Theta$ with $\tau(a)=0$, $\nabla_{\operatorname{ad} a}
= $ left multiplication by $a$.
\end{enumerate}
For simplicity we write $\nabla_j$ for $\nabla_{\partial_j}$.
The operators $\nabla_j\co \mathcal{X}_\Theta \to \mathcal{X}_\Theta$ determine the connection, because of
condition (4), Theorem~\ref{thm:BEJ}, and the linearity axiom, condition~(1).
Once again, the axioms so far correspond to a~connection on the complexif\/ied tangent bundle, so it's
natural to require the covariant derivative of a~``real'' vector f\/ield in a~``real'' direction to be
``real-valued''. In the presence of a~Riemannian metric satisfying Def\/inition~\ref{def:Rmetric}(1)--(4),
this corresponds to the additional axiom
\begin{enumerate}\itemsep=0pt
\setcounter{enumi}{4} \item For any $j$, $k$, and $\ell$, $\langle\nabla_j\partial_k, \partial_\ell
\rangle$ is \emph{self-adjoint}.
\end{enumerate}

We call the connection \emph{torsion-free} if for all $j,k\le n$, $\nabla_j\partial_k = \nabla_k\partial_j$.
This is the exact analogue of the corresponding condition in the commutative case (since $\partial_k$ and
$\partial_j$ commute).

We say the connection is \emph{compatible with a~Riemannian metric} $g=\langle \cdot, \cdot\rangle$ (in the
sense of Def\/inition~\ref{def:Rmetric}) if for all $X,Y\in \mathcal{X}_\Theta$, $Z\in \mathcal{D}_\Theta$,
\begin{gather*}
Z\cdot\langle X,Y\rangle=\langle\nabla_Z X,Y\rangle+\langle X,\nabla_Z Y\rangle.
\end{gather*}
\end{Definition}
\begin{Remark}
Note that a~def\/inition of connections on ``vector bundles'' over noncommutative tori (and other
noncommutative spaces also equipped with a~Lie group of symmetries) was given about 30 years ago by Connes
in his classic paper~\cite{MR572645}.
This def\/inition is also based on~\eqref{eq:conn} (except with left modules replaced by right modules),
but with the vector f\/ield $X$ restricted to have ``constant coef\/f\/icients,'' or in our situation, to
be a~$\mathbb{C}$-linear combination of the $\partial_j$.
In this same paper Connes gives the def\/inition of compatibility with a~metric, and it is the same as ours.
However, he does not address the notion of torsion for a~connection, nor does he attempt to prove a~version
of Levi-Civita's theorem.
\end{Remark}

\section{Levi-Civita's theorem}
Now we can state and prove the analogue of Levi-Civita's theorem.
\begin{Theorem}
\label{thm:LC}
Let $\Theta$ be a~generic skew-symmetric $n\times n$ matrix in the sense of Theorem~\textup{\ref{thm:BEJ}}
and let $g=\langle \cdot, \cdot\rangle$ be a~Riemannian metric on $A^\infty_\Theta$ in the sense of
Definition~\textup{\ref{def:Rmetric}} $($including condition~\textup{\ref{def:Rmetric}(4)}$)$.
Then there is one and only one connection on $\mathcal{X}_\Theta$ in the sense of Definition~\textup{\ref{def:conn}} $($including conditions~\textup{\ref{def:conn}(4)} and~\textup{\ref{def:conn}(5)}$)$ that is torsion-free and compatible with the metric.
This connection, called the {\bfseries Levi-Civita connection}, is determined by the formula
\begin{gather}
\label{eq:LC}
\langle\nabla_j\partial_k,\partial_\ell\rangle
=\frac12\big[\partial_j\langle\partial_k,\partial_\ell\rangle+\partial_k\langle\partial_j,
\partial_\ell\rangle-\partial_\ell\langle\partial_j,\partial_k\rangle\big].
\end{gather}
\end{Theorem}
\begin{proof}
First we prove uniqueness.
Suppose we have a~torsion-free connection $\nabla$ compatible with the metric.
We have (because of compatibility with the metric)
\begin{gather}
\partial_j\langle\partial_k,\partial_\ell\rangle
=\langle\nabla_j\partial_k,\partial_\ell\rangle+\langle\partial_k,\nabla_j\partial_\ell\rangle,
\nonumber
\\
\partial_k\langle\partial_\ell,\partial_j\rangle
=\langle\nabla_k\partial_\ell,\partial_j\rangle+\langle\partial_\ell,\nabla_k\partial_j\rangle,
\label{eq:compatibility}
\\
\partial_\ell\langle\partial_j,\partial_k\rangle
=\langle\nabla_\ell\partial_j,\partial_k\rangle+\langle\partial_j,\nabla_\ell\partial_k\rangle.
\nonumber
\end{gather}
Via conditions~\ref{def:Rmetric}(2),~\ref{def:Rmetric}(4) and~\ref{def:conn}(5), together with the
torsion-free condition, we can re\-wri\-te~\eqref{eq:compatibility} as
\begin{gather}
\label{eq:compatibility1a}
\langle\nabla_j\partial_k,\partial_\ell\rangle
=\partial_j\langle\partial_k,\partial_\ell\rangle-\langle\nabla_\ell\partial_j,\partial_k\rangle,
\\
\label{eq:compatibility1b}
\langle\nabla_j\partial_k,\partial_\ell\rangle
=\partial_k\langle\partial_j,\partial_\ell\rangle-\langle\nabla_\ell\partial_k,\partial_j\rangle,
\qquad
\text{and}
\\
\label{eq:compatibility1c}
0=-\partial_\ell\langle\partial_j,\partial_k\rangle
+\langle\nabla_\ell\partial_j,\partial_k\rangle
+\langle\nabla_\ell\partial_k,\partial_j\rangle.
\end{gather}
Adding~\eqref{eq:compatibility1a},~\eqref{eq:compatibility1b}, and~\eqref{eq:compatibility1c}
gives~\eqref{eq:LC}.

\looseness=-1
Next we prove existence, by showing that~\eqref{eq:LC} determines a~unique connection which is compatible
with the metric and torsion-free.
To begin with, given $j$ and $k$, knowing $\langle \nabla_j\partial_k, \partial_\ell \rangle$ for all
$\ell$ determines $\nabla_j\partial_k$, since the metric is nondegenerate.
We get a~unique extension to a~def\/inition of $\nabla_X\partial_k$ for all vector f\/ields
$X\in\mathcal{D}_\Theta$ by making $\nabla_X\partial_k$ linear in $X$ and requiring that
$\nabla_{\operatorname{ad} a}\partial_k = a\partial_k$ when $\tau(a)=0$.
(We are using genericity of $\Theta$ in order to appeal to Theorem~\ref{thm:BEJ}.) Then since the
$\partial_k$ are a~free $A_\Theta^\infty$-basis for $\mathcal{X}_\Theta$, knowing $\nabla_X\partial_k\in
\mathcal{X}_\Theta$ for each $k$ uniquely determines $\nabla_X Y\in \mathcal{X}_\Theta$ for each $Y\in
\mathcal{X}_\Theta$, since we have a~unique expression $Y = \sum\limits_k a_k \partial_k$ and can set
\begin{gather*}
\nabla_X\left(\sum_k a_k\partial_k\right)=\sum_k(X\cdot a_k)\partial_k+a_k\nabla_X\partial_k.
\end{gather*}
This gives us a~def\/inition of $\nabla$ satisfying the axioms of Def\/inition~\ref{def:conn}(1),~(2)
and~\ref{def:conn}(4).
Condition~\ref{def:conn}(5) holds because of condition~\ref{def:Rmetric}(4) and the fact that the
$\partial_j$ are $*$-preserving.
We need to show that~\ref{def:conn}(3) is also satisf\/ied, which means we need to check that
\begin{gather}
\nabla_X\bigl(ab\partial_k\bigr)=(X\cdot a)\bigl(b\partial_k\bigr)+a\nabla_X\bigl(b\partial_k\bigr).
\label{eq:textconn}
\end{gather}
The left-hand side of~\eqref{eq:textconn} is def\/ined to be
\begin{gather*}
\bigl(X\cdot(ab)\bigr)\partial_k+ab\nabla_X\partial_k.
\end{gather*}
Since $X\in \mathcal{D}_\Theta$ and is thus a~derivation, this becomes
\begin{gather*}
(X\cdot a)\bigl(b\partial_k\bigr)+a(X\cdot b)\partial_k+ab\nabla_X\partial_k
=(X\cdot a)\bigl(b\partial_k\bigr)+a\nabla_X\bigl(b\partial_k\bigr),
\end{gather*}
which agrees with the right-hand side of~\eqref{eq:textconn}, as required.

The right-hand side of~\eqref{eq:LC} is clearly symmetric under interchange of $j$ and $k$, because of the
fact that Def\/inition~\ref{def:Rmetric}(4) ensures that $\langle \partial_j, \partial_k\rangle = \langle
\partial_k, \partial_j\rangle$, so $\nabla$ is torsion-free.
We have just one more thing to check, which is compatibility with the metric.
From~\eqref{eq:LC}, we have
\begin{gather*}
\langle\nabla_j\partial_k,\partial_\ell\rangle+\langle\partial_k,\nabla_j\partial_\ell\rangle
=\frac12\Bigl[\partial_j\langle\partial_k,
\partial_\ell\rangle+\partial_k\langle\partial_j,\partial_\ell\rangle-\partial_\ell\langle\partial_j,\partial_k\rangle
\\
\phantom{\langle\nabla_j\partial_k,\partial_\ell\rangle+\langle\partial_k,\nabla_j\partial_\ell\rangle=}
{}+\partial_j\langle\partial_\ell,\partial_k\rangle+\partial_\ell\langle\partial_j,
\partial_k\rangle-\partial_k\langle\partial_j,\partial_\ell\rangle\Bigr]
=\partial_j\langle\partial_k,\partial_\ell\rangle,
\end{gather*}
which is what is required.
This completes the proof.
\end{proof}

\section{Riemannian curvature}

With Theorem~\ref{thm:LC} in place, it now makes sense to def\/ine curvature for a~Riemannian metric using
derivatives of the Levi-Civita connection.
There are dif\/ferent sign conventions used by dif\/ferent authors; here we are following~\cite{MR1138207}.
\begin{Definition}
\label{def:curv}
Let $\Theta$ be a~generic skew-symmetric $n\times n$ matrix in the sense of Theorem~\textup{\ref{thm:BEJ}}
and let $g=\langle \cdot, \cdot\rangle$ be a~Riemannian metric on $A^\infty_\Theta$ in the sense of
Def\/inition~\textup{\ref{def:Rmetric}} $($including condition~\textup{\ref{def:Rmetric}(4)}$)$.
Let $\nabla$ be the associated Levi-Civita connection from Theorem~\ref{thm:LC}.
Def\/ine the associated \emph{Riemann curvature operator} by
\begin{gather*}
R(X,Y)=\nabla_Y\nabla_X-\nabla_X\nabla_Y+\nabla_{[X,Y]}\co\mathcal{X}_\Theta\to\mathcal{X}_\Theta,
\qquad
X,Y\in\mathcal{D}_\Theta.
\end{gather*}
We also def\/ine the \emph{Riemannian curvature} by
\begin{gather*}
R_{j,k,\ell,m}=\left\langle R(\partial_j,\partial_k)\partial_\ell,\partial_m\right\rangle.
\end{gather*}
\end{Definition}
\begin{Remark}
Recall by Theorem~\ref{thm:BEJ} that $\mathcal{D}_\Theta$ splits as the direct sum of the $n$-dimensional
vector space spanned by the $\partial_j$ and the set of $\operatorname{ad} a$, $a\in A^\infty_\Theta$.
The second summand actually has no ef\/fect on the curvature the way we've normalized the connection, for
if $a\in A^\infty_\Theta$, $\tau(a)=0$, and if $X\in \mathcal{D}_\Theta$, then f\/irst of all $\tau(Xa)=0$
(this is obvious for $X$ inner, so we only need to check it for $X=\partial_j$, where it follows from the
fact that the gauge action of $\mathbb{T}^n$ on $A_\Theta$ preserves $\tau$~--- see also~\cite[Lemma
2.1]{MR2444094}).
So for $b\in A^\infty_\Theta$,
\begin{gather*}
\lbrack\operatorname{ad}(a),X\rbrack b=[a,X\cdot b]-X\cdot([a,b])
=a(X\cdot b)-(X\cdot b)a-X\cdot(ab-ba)
\\
\phantom{\lbrack\operatorname{ad}(a),X\rbrack b}
=a(X\cdot b)-(X\cdot b)a-(X\cdot a)b-a(X\cdot b)+(X\cdot b)a-b(X\cdot a)
\\
\phantom{\lbrack\operatorname{ad}(a),X\rbrack b}
=-[X\cdot a,b]=\operatorname{ad}(-X\cdot a)b,
\end{gather*}
and we have $[\operatorname{ad} a, X]=\operatorname{ad}(-X\cdot a)$.
So
\begin{gather*}
R(\operatorname{ad}a,X)Z=\left(\nabla_X\nabla_{\operatorname{ad}a}-\nabla_{\operatorname{ad}a}
\nabla_X+\nabla_{[\operatorname{ad}a,X]}\right)Z
\\
\phantom{R(\operatorname{ad}a,X)Z}
=\nabla_X(a Z)-a(\nabla_X Z)+\nabla_{\operatorname{ad}(-X\cdot a)}Z
\\
\phantom{R(\operatorname{ad}a,X)Z}
=(X\cdot a)Z+a(\nabla_X Z)-a(\nabla_X Z)+(-X\cdot a)Z=0.
\end{gather*}
Since $R(X,Y)$ is bilinear and antisymmetric in $X$ and $Y$, it follows that $R(X,Y)$ only depends on the
projections of $X$ and $Y$ into the $\mathbb{C}$-span of $\partial_1,\dots,\partial_n$.
\end{Remark}
\begin{Proposition}
In the context of Definition \textup{\ref{def:curv}}, if $X,Y \in \mathcal{D}_\Theta$, then $R(X,Y)$ is
$A^\infty_\Theta$-linear, i.e., ``is a~tensor''.
\end{Proposition}
\begin{proof}
The classical proof works without change.
Just expand
\begin{gather*}
R(X,Y)(aZ)=\bigl(\nabla_Y\nabla_X-\nabla_X\nabla_Y+\nabla_{[X,Y]}\bigr)\bigl(aZ\bigr)
\end{gather*}
using the Leibniz rule and observe that all the terms involving $X\cdot a$ and $Y\cdot a$ cancel out.
\end{proof}
\begin{Proposition}\label{prop:curvprops}
In the context of Definition \textup{\ref{def:curv}}, the Riemannian curvature satisfies the following
symmetry properties for all $j$, $k$, $\ell$, $m$:
\begin{enumerate}\itemsep=0pt
\item[$1.$] $R_{j,k,\ell,m} + R_{k,\ell,j,m} + R_{\ell,j,k,m} = 0$ $(${\bfseries Bianchi
identity}$)$;
\item[$2.$] $R_{j, k, \ell,m} = - R_{k, j, \ell,m}$.
\end{enumerate}
\end{Proposition}
\begin{proof}
The easiest is (2), which is immediate from the fact that the def\/inition of $R(X, Y)$ is antisymmetric in
$X$ and $Y$.

Next we prove (1).
We expand $R(\partial_j, \partial_k)\partial_\ell$, etc., using the fact that $\nabla$ is torsion-free, and
obtain
\begin{gather*}
R(\partial_j,\partial_k)\partial_\ell+R(\partial_k,\partial_\ell)\partial_j+R(\partial_\ell,\partial_j)\partial_k
\\
\qquad
=\nabla_k\nabla_j\partial_\ell-\nabla_j\nabla_k\partial_\ell
+\nabla_\ell\nabla_k\partial_j-\nabla_k\nabla_\ell\partial_j
+\nabla_j\nabla_\ell\partial_k-\nabla_\ell\nabla_j\partial_k
\\
\qquad
=(\nabla_k\nabla_\ell\partial_j-\nabla_k\nabla_\ell\partial_j)
+(\nabla_\ell\nabla_j\partial_k-\nabla_\ell\nabla_j\partial_k)
+(\nabla_j\nabla_k\partial_\ell-\nabla_j\nabla_k\partial_\ell)
=0,
\end{gather*}
proving the Bianchi identity.
\end{proof}

As in classical Riemannian geometry, the symmetry properties of the curvature
(Proposition~\ref{prop:curvprops}) greatly cut down the number of
independent components of the curvature, especially in low dimension.
But two additional symmetry properties that hold in the commutative case,
$R_{j,k,\ell,m}=-R_{j,k,m,\ell}$ and $R_{j,k,\ell,m}=R_{\ell,m,j,k}$,
fail in general, as pointed out to me by Joakim Arnlind.
For example, in the example below
for the irrational rotation algebra~$A_\theta$ (the case~$n=2$),
$R_{1,2,2,2}$ and $R_{1,2,1,1}$ can be nonzero.

\section{An example}
We illustrate our theory in the case of the irrational rotation algebra $A_\theta$ (for generic $\theta$),
and a~Riemannian metric that is a~conformal deformation of the f\/lat metric associated to the complex
elliptic curve $\mathbb{C}/(\mathbb{Z}+i\mathbb{Z})$.
For this f\/lat metric, we have $\langle \partial_j, \partial_k\rangle = \delta_{jk}$, so we choose
a~conformal factor $h=h^*\in A_\theta^\infty$ and suppose $\langle \partial_j, \partial_k\rangle =
e^h\delta_{jk}$.
Formula~\eqref{eq:LC} then determines the Levi-Civita connection; we have
\begin{gather}
\label{eq:T2conn}
\langle\nabla_j\partial_k,\partial_l\rangle=\frac12\bigl\lbrack\delta_{kl}\partial_j\big(e^h\big)+\delta_{jl}
\partial_k\big(e^h\big)-\delta_{jk}\partial_\ell\big(e^h\big)\bigr\rbrack.
\end{gather}
For example, from~\eqref{eq:T2conn},
\begin{gather*}
\langle\nabla_1\partial_1,\partial_1\rangle=\frac12\partial_1\big(e^h\big),
\qquad
\langle\nabla_2\partial_2,\partial_2\rangle=\frac12\partial_2\big(e^h\big).
\end{gather*}
We also get
\begin{gather*}
\langle\nabla_1\partial_1,\partial_2\rangle=-\frac12\partial_2\big(e^h\big),
\qquad
\langle\nabla_2\partial_2,\partial_1\rangle=-\frac12\partial_1\big(e^h\big).
\end{gather*}
So these imply that
\begin{gather*}
\nabla_1\partial_1=-\nabla_2\partial_2=\frac12\bigl(\partial_1\big(e^h\big)e^{-h}\partial_1-\partial_2\big(e^h\big)e^{-h}
\partial_2\bigr)=\frac12\bigl(k_1\partial_1-k_2\partial_2\bigr),
\end{gather*}
where we have written $k_j = \partial_j(e^h)e^{-h}$.
Similarly
\begin{gather*}
\nabla_2\partial_1=\nabla_1\partial_2=\frac12\bigl(\partial_2\big(e^h\big)e^{-h}\partial_1+\partial_1\big(e^h\big)e^{-h}
\partial_2\bigr)=\frac12\bigl(k_2\partial_1+k_1\partial_2\bigr).
\end{gather*}
This makes it possible to compute the curvature.
We obtain
\begin{gather*}
R_{1,2,1,2}=\langle R(\partial_1,\partial_2)\partial_1,\partial_2\rangle
=\langle\nabla_2\nabla_1\partial_1-\nabla_1\nabla_2\partial_1,\partial_2\rangle
\\
\phantom{R_{1,2,1,2}}
=\frac12\Bigl\langle\nabla_2\bigl(k_1\partial_1-k_2\partial_2\bigr)
-\nabla_1\bigl(k_2\partial_1+k_1\partial_2\bigr),\partial_2\Bigr\rangle
\\
\phantom{R_{1,2,1,2}}
=\frac12\Bigl\langle\partial_2(k_1)\partial_1+k_1\nabla_2\partial_1-\partial_2(k_2)\partial_2-k_2\nabla_2\partial_2
\\
\phantom{R_{1,2,1,2}=}
{}-\partial_1(k_2)\partial_1-k_2\nabla_1\partial_1-\partial_1(k_1)\partial_2-k_1\nabla_1\partial_2,\partial_2\Bigr\rangle.
\end{gather*}
The four terms without a~$\nabla_j$ in them contribute
\begin{gather*}
-\frac12\bigl\langle\bigl(\partial_2(k_2)+\partial_1(k_1)\bigr)\partial_2,\partial_2\bigr\rangle
=-\frac12\bigl(\partial_2(k_2)+\partial_1(k_1)\bigr)e^h.
\end{gather*}
The remaining four terms contribute
\begin{gather*}
\frac12\Bigl[k_1\langle\nabla_2\partial_1,\partial_2\rangle-k_2\langle\nabla_2\partial_2,\partial_2\rangle
-k_2\langle\nabla_1\partial_1,\partial_2\rangle-k_1\langle\nabla_1\partial_2,\partial_2\rangle\Bigr]
\\
\qquad
=\frac14\Bigl[k_1\partial_1\big(e^h\big)-k_2\partial_2\big(e^h\big)+k_2\partial_2\big(e^h\big)-k_1\partial_1\big(e^h\big)\Bigr]=0.
\end{gather*}
So we conclude that
\begin{gather}
R_{1,2,1,2}=-\frac12\bigl(\partial_2(k_2)+\partial_1(k_1)\bigr)e^h
\label{eq:R1212a}
\end{gather}
and expanding using the definitions of $k_1$ and $k_2$:
\begin{gather}
\phantom{R_{1,2,1,2}}
=-\frac12\bigl(\partial_2^2\big(e^h\big)e^{-h}-\partial_2\big(e^h\big)e^{-h}\partial_2\big(e^h\big)e^{-h}\notag
+\partial_1^2\big(e^h\big)e^{-h}-\partial_1\big(e^h\big)e^{-h}\partial_1\big(e^h\big)e^{-h}\big)e^h\notag
\nonumber
\\
\phantom{R_{1,2,1,2}}
=-\frac12\big(\Delta\big(e^h\big)-\partial_1\big(e^h\big)e^{-h}\partial_1\big(e^h\big)-\partial_2\big(e^h\big)e^{-h}\partial_2\big(e^h\big)\big),
\end{gather}
where $\Delta$ is the Laplacian $\partial_1^2 + \partial_2^2$.
If $h$ and its derivatives all commute, then $k_j$ would just be~$\partial_j(h)$ and this would reduce to
$-\frac12 e^h\Delta h$.
On the other hand, in the commutative case, we would really want the \emph{Gaussian curvature}~$K$, which
would be $e^{-2h}R_{1,2,1,2}$ (since the vector f\/ields~$\partial_1$ and~$\partial_2$ are orthogonal but
not normalized), and we'd get the classical formula $K=-\frac12 e^{-h}\Delta h$.
Our calculation is clearly related to, but vastly simpler, than the calculations
in~\cite{2011arXiv1110.3500C} and~\cite{2011arXiv1110.3511F}.
Reconciling these very dif\/ferent approaches to curvature in noncommutative geometry is an important
problem for the future.
However, we note that we do have an analogue of the Gauss--Bonnet theorem in our context, which can be
formulated as follows:
\begin{Proposition}[Gauss--Bonnet theorem]
\label{prop:GB}
Let $A_\theta^\infty$ be a~smooth irrational rotation algebra, with generic $\theta$, equipped with
a~Riemannian metric $\langle \partial_j, \partial_k\rangle = e^h\delta_{jk}$, $h=h^*\in A_\theta^\infty$ as
above.
Then if $\tau $ is the canonical trace on $A_\theta$, we have $\tau(R_{1,2,1,2} e^{-h})=0$.
\end{Proposition}
\begin{proof}
By formula~\eqref{eq:R1212a}, $R_{1,2,1,2}\,e^{-h}$ is, up to a~factor of $-\frac12$, just
$\partial_1(k_1)+\partial_2(k_2)$.
But $\tau(\partial_j(a)) = 0$ for any $a$, since $\tau$ is invariant under the gauge action of
$\mathbb{T}^2$ (see also~\cite[Lemma 2.1]{MR2444094}).
\end{proof}
\begin{Remark}
We should explain why Proposition~\ref{prop:GB}, in the commutative case of $\mathbb{T}^2$ ($\theta=0$),
really is the Gauss--Bonnet theorem.
In that case, $\tau$ is integration against Haar measure, and so $\tau(\,\cdot\,e^h)$ is integration
against the Riemannian volume form, which dif\/fers from the standard volume form by $\sqrt{\det(g)} = e^h$.
Since $K=e^{-2h}R_{1,2,1,2}$, the integral of $K$ against the Riemannian volume form is thus
$\tau(R_{1,2,1,2}\,e^{-2h}e^h) = \tau(R_{1,2,1,2}\,e^{-h})$.
\end{Remark}

\subsection*{Acknowledgements} This research was supported by NSF grant DMS-1206159.
The author thanks the referees and the participants at the Fields Institute Focus Program for several
interesting comments and discussions. I~would like to thank Joakim Arnlind for pointing out a mistake in the
original formulation of Proposition~\ref{prop:curvprops}.

\pdfbookmark[1]{References}{ref}
\LastPageEnding

\end{document}